\documentclass[12pt]{article}
\usepackage{amssymb,amsmath, amsthm,latexsym, verbatim, amscd}
\usepackage{here}

  \newtheorem{theorem}{Theorem}[section] %

  \newtheorem{definition-theorem}[theorem]{Definition-Theorem}
\theoremstyle{definition} %
  \newtheorem{definition}[theorem]{Definition} %
  \newtheorem{example}[theorem]{Example} %
  \newtheorem{problem}[theorem]{Problem}

  \newtheorem{conjecture}[theorem]{Conjecture}
  \newtheorem{remark}[theorem]{Remark} %
  \newtheorem{setting}[theorem]{Setting}

\newcommand{\invHom}[3]{\operatorname{Hom}_{#1}(#2,#3)}

\begin{document}
\title{Conjectures on reductive homogeneous spaces}
\author{Toshiyuki KOBAYASHI
\footnote{Graduate School of Mathematical Sciences, 
The University of Tokyo, 
3-8-1 Komaba, Tokyo 153-8914, Japan.  }
\date{}
}
\maketitle

\begin{abstract}
We address some conjectures
 and open problems 
 in \lq{analysis of symmetries}\rq\
which include the study of non-commutative harmonic analysis 
 and discontinuous groups
 for reductive homogeneous spaces
 beyond the classical framework:
\newline
(1)
discrete series for non-symmetric homogeneous spaces $G/H$;
\newline
(2)
discontinuous group $\Gamma$ for $G/H$
 beyond the Riemannian setting;
\newline
(3)
analysis on pseudo-Riemannian locally homogeneous spaces.  
\end{abstract}


\setcounter{section}{0}
\section{Introduction}

I have been working on various subjects of mathematics,
 and \lq\lq{symmetry}\rq\rq\
is a key word to create new interactions among different disciplines.

\vskip 1pc

In this  paper, 
 I would like to address some conjectures and problems
 of the areas of  \lq\lq{analysis of symmetries}\rq\rq\
 in which I have been deeply involved, 
 but to which I do not have a solution.   

\begin{enumerate}
\item[(1)]
Discrete series for non-symmetric homogeneous spaces $G/H$;
\item[(2)]
Discontinuous group $\Gamma$ for $G/H$
 beyond the Riemannian setting;
\item[(3)]
Analysis on non-Riemannian locally homogeneous spaces $\Gamma \backslash G/H$.  
\end{enumerate}

These three topics are discussed
 in Sections \ref{sec:discrep}--\ref{sec:spec}, 
 respectively, 
 using a common setting
 which is explained
 in Section \ref{sec:setting}
 with simple examples.

\section{Basic settings}
\label{sec:setting}
Throughout this paper, 
 our basic geometric setting will be as follows.  
\begin{setting}
\label{set:reductive}
$G$ is a real reductive linear Lie group, 
 $H$ is a closed proper subgroup
 which is reductive in $G$, 
 and $X:=G/H$.  
\end{setting}

A distinguished feature
 of this setting is that the manifold $X$ carries a {\it{pseudo-Riemannian structure}} 
 with a \lq{large}\rq\ isometry group, 
 namely, 
 the reductive group $G$ acts transitively and isometrically
 on $X$.  
Such a pseudo-Riemannian structure is induced from 
 the Killing form $B$
 if $G$ is semisimple.  
For reductive $G$, 
 one can take a maximal compact subgroup $K$ of $G$
 such that $H \cap K$ is a maximal compact subgroup of $H$, 
 and a $G$-invariant symmetric bilinear form $B$
 on the Lie algebra ${\mathfrak{g}}$ of $G$
 such that the Cartan decomposition 
 ${\mathfrak{g}}={\mathfrak{k}}+{\mathfrak{p}}$
 is an orthogonal decomposition with respect to $B$
 and that $B$ is negative definite 
 on ${\mathfrak{k}}$ 
 and is positive definite 
 on ${\mathfrak{p}}$.  
Then $B$ induces a $G$-invariant pseudo-Riemannian structure of signature $(p,q)$
 on $X$, 
 where $p+q= \dim X$
 and $q=\dim K/H\cap K$.

Very special cases
 of homogeneous spaces 
 in Setting \ref{set:reductive} include:

\begin{example}
[semisimple coadjoint orbits]
For a reductive Lie group $G$, 
 one can identify the Lie algebra ${\mathfrak{g}}$ 
 with its dual ${\mathfrak{g}}^{\ast}$ via $B$.  
The coadjoint orbit ${\mathcal{O}}_{\lambda}:=\operatorname{Ad}^{\ast}(G) \lambda$
 is called {\it{semisimple}}, 
 {\it{elliptic}}, 
 {\it{hyperbolic}}, 
 or {\it{regular}}
 if the element in ${\mathfrak{g}}$
 corresponding to $\lambda$ has that property.  
We write ${\mathfrak{g}}_{\operatorname{ss}}^{\ast}$,
 ${\mathfrak{g}}_{\operatorname{ell}}^{\ast}$,
 ${\mathfrak{g}}_{\operatorname{hyp}}^{\ast}$,
 or ${\mathfrak{g}}_{\operatorname{reg}}^{\ast}$
 for the collection of such elements $\lambda$, 
 respectively.  
By definition,
 ${\mathfrak{g}}_{\operatorname{ell}}^{\ast}$,
 ${\mathfrak{g}}_{\operatorname{hyp}}^{\ast}$
 $\subset {\mathfrak{g}}_{\operatorname{ss}}^{\ast}$.  
The isotropy subgroup of $\lambda$ 
 is reductive
 if $\lambda \in {\mathfrak{g}}_{\operatorname{ss}}^{\ast}$, 
 hence any semisimple coadjoint orbit ${\mathcal{O}}_{\lambda}$
 gives an example of Setting \ref{set:reductive}.

For compact $G$, 
 one has 
$
   {\mathfrak{g}}_{\operatorname{ell}}^{\ast}
  ={\mathfrak{g}}_{\operatorname{ss}}^{\ast}
  ={\mathfrak{g}}^{\ast}
$
 and ${\mathcal{O}}_{\lambda}$ is a generalized flag variety
 for any $\lambda \in {\mathfrak{g}}^{\ast}$;
 ${\mathcal{O}}_{\lambda}$ is a full flag variety
 iff $\lambda \in {\mathfrak{g}}_{\operatorname{reg}}^{\ast}$.

The subclass $\{{\mathcal{O}}_{\lambda} : \lambda \in {\mathfrak{g}}_{\operatorname{ss}}^{\ast}\}$ in Setting \ref{set:reductive} plays a particular role
 in the unitary representation theory.  
The orbit philosophy
 due to Kirillov--Kostant--Duflo--Vogan suggests 
 an intimate relationship 
 between the set ${\mathfrak{g}}^{\ast}/\operatorname{Ad}^{\ast}(G)$
 of coadjoint orbits
 and the set of equivalence classes 
 of irreducible unitary representations of $G$
 (the {\it{unitary dual}} $\widehat G$):

\begin{equation}
\label{eqn:orbit}
{\mathfrak{g}}^{\ast}/\operatorname{Ad}^{\ast}(G)
\fallingdotseq\widehat G, 
\qquad
{\mathcal{O}}_{\lambda} \leftrightarrow \pi_{\lambda}.  
\end{equation}
We recall that any coadjoint orbit carries a canonical symplectic form
 called the {\it{Kirillov--Kostant--Souriau form}}.  
Then the correspondence ${\mathcal{O}}_{\lambda} \mapsto \pi_{\lambda}$
 is supposed to be a \lq\lq{geometric quantization}\rq\rq\
 of the Hamiltonian $G$-manifold ${\mathcal{O}}_{\lambda}$
 if such $\pi_{\lambda}$ exists.  
This philosophy works fairly well 
 for $\lambda \in {\mathfrak{g}}_{\operatorname{ss}}^{\ast}$
 satisfying an appropriate integral condition:
 loosely speaking, 
 $\pi_{\lambda}$ is obtained by a unitary induction from 
 a parabolic subgroup 
 (real polarization of the para-Hermitian manifold ${\mathcal{O}}_{\lambda}$)
 for $\lambda \in {\mathfrak{g}}_{\operatorname{hyp}}^{\ast}$, 
 in a Dolbeault cohomology space
 on the pseudo-K{\"a}hler manifold ${\mathcal{O}}_{\lambda}$
 as a generalization
 of the Borel--Weil--Bott theorem 
 (complex polarization of ${\mathcal{O}}_{\lambda}$)
 or alternatively
 by a cohomological parabolic induction
 ({\it{e.g., }} Zuckerman's derived functor module $A_{\mathfrak{q}}(\lambda)$)
 for $\lambda \in {\mathfrak{g}}_{\operatorname{ell}}^{\ast}$, 
 and by the combination
 of these two procedures
 for general $\lambda \in {\mathfrak{g}}_{\operatorname{ss}}^{\ast}$,
 although there are some delicate issues
 about singular $\lambda$
 and also about \lq\lq{$\rho$-shift}\rq\rq\
 of the parameter, 
 see \cite[Chap.\ 2]{xkrons94}
 for instance.  
The resulting 
 \lq\lq{quantizations}\rq\rq\
 $\pi_{\lambda}$
 of semisimple coadjoint orbits ${\mathcal{O}}_{\lambda}$ 
 give \lq\lq{large part}\rq\rq\ of the unitary dual $\widehat G$.  
\end{example}

\begin{example}
[symmetric spaces, real spherical spaces]
\label{ex:symm}
Let $\sigma$ be an automorphism
 of a reductive Lie group $G$
 of finite order, 
 $G^{\sigma}$ the fixed point subgroup 
 of $\sigma$, 
 and $H$ an open subgroup of $G^{\sigma}$.  
Then $H$ is reductive
 and the homogeneous space $G/H$ provides another example
 of Setting \ref{set:reductive}.  
In particular, 
 if the order of $\sigma$ is two, 
 $G/H$ is called a (reductive) {\it{symmetric space}}.  
Geometrically,
 it is a symmetric space 
 with respect to the Levi-Civita connection 
 of the pseudo-Riemannian structure
 in the sense 
 that all geodesic symmetries are globally defined isometries.  
This is a subclass of Setting \ref{set:reductive}
 for which the $L^2$-analysis has been
 extensively studied over 60 years.  
Group manifolds
 $(G \times G)/\operatorname{diag}(G)$, 
 Riemannian symmetric spaces $G/K$
 and irreducible affine symmetric spaces
 such as $SL(p+q, {\mathbb{R}})/SO(p,q)$
 are examples
 of reductive symmetric spaces.  
More generally, 
 in Setting \ref{set:reductive},
 one has
\[
  \{\text{symmetric spaces}\}
  \subset
  \{\text{spherical spaces}\}
  \subset
  \{\text{real spherical spaces}\}, 
\]
where we say $G/H$ is {\it{spherical}}
 if a Borel subgroup of the complexification $G_{\mathbb{C}}$ has an open orbit
 in $G_{\mathbb{C}}/H_{\mathbb{C}}$, 
 and $G/H$ is {\it{real spherical}}
 if a minimal parabolic subgroup of $G$
 has an open orbit
 in $G/H$. 
See Kobayashi--T.\ Oshima  \cite{KOt13}
for the roles that these geometric properties play in the global analysis on $G/H$.  
When $H$ is compact, 
 $G/H$ is spherical
 if and only if it is a weakly symmetric space
 in the sense of Selberg.  
\end{example}

The model space of non-zero constant sectional curvatures
 in pseudo-Riemannian geometry
 is a special case
 of reductive symmetric spaces:

\begin{example}
[pseudo-Riemannian space form, 
 see {\cite{W11}}]
\label{ex:1.4}
The hypersurface
\[
   X(p,q):=\{x \in {\mathbb{R}}^{p+q+1}
   :x_1^2+\cdots+x_{p+1}^2-x_{p+2}^2-\cdots-x_{p+q+1}^2=1\}
\]
in ${\mathbb{R}}^{p+1,q}:=({\mathbb{R}}^{p+q+1}, 
 ds^2 = dx_1^2+\cdots + dx_{p+1}^2-d x_{p+2}^2-\cdots-dx_{p+q}^2)$
 carries a pseudo-Riemannian structure
 of signature $(p,q)$ with constant sectional curvature 1.  
Equivalently,
 we may regard $X(p,q)$
 as a space of constant sectional curvature $-1$
 with respect to the pseudo-Riemannian metric of signature $(q,p)$.  
If $q=0$, $p=0$, $q=1$,  or $p=1$, 
 then $X(p,q)$ is the sphere ${\operatorname{S}}^p$, 
 the hyperbolic space ${\operatorname{H}}^q$, 
 the de Sitter space ${\operatorname{d S}}^{p+1}$, 
 or the anti-de Sitter space $\operatorname{A d S}^{q+1}$, 
 respectively.  
For general $(p,q)$, 
 the generalized Lorentz group $O(p+1,q)$ acts transitively
 and isometrically on $X(p,q)$, 
 and one has a diffeomorphism $X(p,q)\simeq O(p+1,q)/O(p,q)$, 
 giving an expression of $X(p,q)$ as a reductive symmetric space of rank one.  
\end{example}

\section{Problems on discrete series for $G/H$}
\label{sec:discrep}

The \lq{smallest units of symmetries}\rq\
 defined by group actions
 may be 
\newline\indent\indent
{\it{irreducible representations}}
 if the action is linear, 
 and 
\newline\indent\indent
{\it{homogeneous spaces}}
 if the action is smooth 
 on a manifold.

The objects of this section are irreducible subrepresentations
 in $L^2(X)$
 for homogeneous spaces $X$, 
 that is, 
 {\it{discrete series representations}}
 for $X$
 (Definition \ref{def:D1} below), 
 which are supposed to be a building block in global analysis on $X$.  
For instance, 
 when $X$ is a reductive symmetric space, 
 parabolic inductions 
of discrete series representations
 for subsymmetric spaces
 yield the full spectrum
 in the Plancherel formula of $L^2(X)$, 
 see \cite{xdelorme} for instance.

This section elucidates the following problem
 in the generality of Setting \ref{set:reductive}
 by using simple examples, 
 and addresses some related conjectures.

\begin{problem}
\label{prob:D1}
Find all discrete series representations
 for $G/H$.  
\end{problem}

Let us fix some notation.  
Suppose a Lie group $G$ acts continuously 
 on a manifold $X$.  
Then one has a natural unitary representation
 ({\it{regular representation}})
 of $G$ on the Hilbert space $L^2(X)$ 
 of $L^2$-sections for the half-density bundle 
 ${\mathcal{L}} := (\wedge {}^{\dim X} T^*X \otimes {or}_X)^\frac{1}{2}$ of $X$
where ${or}_X$ stands for the orientation bundle.

\begin{definition}
\label{def:D1}
An irreducible unitary representation $\pi$ of $G$
 is said to be a {\it{discrete series representation}} for $X$
 if there exists a non-zero continuous $G$-homomorphism from 
 $\pi$ to the regular representation on $L^2(X)$.  
In other words, 
 discrete series representations for $X$ 
 are irreducible subrepresentations 
 realized in closed subspaces of the Hilbert space $L^2(X)$.  
\end{definition}

We denote by $\operatorname{Disc}(X)$
 the set of discrete series representations for $X$.  
It is a (possibly, empty) subset 
 of the unitary dual $\widehat G$
 of the group $G$.

Hereafter, 
 suppose we are in Setting \ref{set:reductive}.  
Then there is a $G$-invariant Radon measure $\mu$
 on the homogeneous space $X=G/H$, 
 hence ${\mathcal{L}}$ is trivial as a $G$-equivariant bundle
 and $L^2(X)$ may be identified with $L^2(X, d \mu)$.

If $G/H$ is spherical
 (see Example \ref{ex:symm}), 
 then the ring 
\[
  {\mathbb{D}}(G/H):=
  \{\text{$G$-invariant differential operators on $G/H$}\}
\]
 is commutative and the multiplicity
 of irreducible representations $\pi$ of $G$
 in the regular representation
 on $C^{\infty}(G/H)$
 is uniformly bounded, 
 and vice versa \cite{KOt13}.  
In this case, 
 the disintegration of the regular representation
 $L^2(X)$
 into irreducibles
 (the Plancherel-type theorem)
 is essentially equivalent 
 to the joint spectral decomposition 
 for the commutative ring ${\mathbb{D}}(G/H)$, 
 and Problem \ref{prob:D1}
 highlights point spectra in $L^2(G/H)$.

Classical examples
 trace back to Gelfand--Graev (1962), Shintani (1967), 
 Molchanov (1968), 
 J.\ Faraut (1979), 
 R.\ S.\ Strichartz (1983)
 and some others
 on the analysis of the space form $X(p,q)$, 
 which we review with some modern viewpoints, 
 see \cite[Thm.\ 2.1]{KAdv21} and references therein.

\begin{example}
\label{ex:Xpqdisc}
Let  $(G,H)=(O(p+1,q),O(p,q))$ 
 and $X=G/H$ as in Example \ref{ex:1.4}.  
Then the ring ${\mathbb{D}}(G/H)$
 is generated by the Laplacian $\Delta_X$, 
 which is not an elliptic differential operator if $p,q>0$.  
We set
\[
  L^2(X)_{\lambda}
  :=
\{f \in L^2(X): \text{$\Delta_X f = \lambda f$ in the weak sense}
\}.  
\]
Then $L^2(X)_{\lambda}$ is a closed subspace
 in $L^2(X)$, 
 and the resulting unitary representation of $G$
 on $L^2(X)_{\lambda}$ is irreducible
 whenever it is non-zero.  
Conversely, 
 any discrete series representation 
 for $X$ is realized on an $L^2$-eigenspace $L^2(X)_{\lambda}$
 for some eigenvalue $\lambda$.  
In particular, 
 one has the equivalence:
\begin{alignat*}{2}
\operatorname{Disc}(G/H) = \emptyset
&\iff
\text{there is no point spectrum of $\Delta_X$ in $L^2(X)$}
\\
&\iff 
\text{$p=0$ and $q \ge 1$.}
\end{alignat*}
Thus there exists point spectrum
 of the Laplacian $\Delta_X$ in $L^2(X)$
 unless $X = X(p,q)$ is a hyperbolic space
 $\operatorname{H}^q \equiv X(0,q)$.  
The description of the eigenspace $L^2(X)_{\lambda}$
 for $q=0$ is the classical theory 
 of spherical harmonics on the sphere $\operatorname{S}^p \equiv X(p,0)$.  
For $p,q \ge 1$, 
 one has 
\[
  \text{$L^2(X)_{\lambda} \ne 0$
 iff
 $\lambda = \lambda_k$
 for some $k \in {\mathbb{Z}}$
 with $-\frac 1 2 (p+q-1)<k$,  }
\]
 where $\lambda_k:=-k(k+p+q-1)$.  
The resulting irreducible unitary representation
 on $L^2(X)_{\lambda_k}$
 is isomorphic to a \lq{geometric quantization}\rq\ of an elliptic coadjoint orbit
 of minimal dimension, 
 or alternatively 
 in an algebraic language, 
 it is the unitarization
 of Zuckerman's derived functor module
 $A_{\mathfrak{q}}(k)$
 with the normalization as in \cite{VZ84}.  
This algebraic description involves
 delicate questions for finitely many exceptional parameters,
 {\it{i.e., }}
 those for $k<0$, 
 see Problem \ref{prob:D3} below.  
\end{example}

In the generality of Setting \ref{set:reductive}, 
 we may divide Problem \ref{prob:D1} 
 into two subproblems:
\par\noindent
(A) a characterization of the pairs $(G, H)$
 for which $G/H$ admits at least one discrete series representation
 (Problem \ref{prob:D2});
\par\noindent
(B) a description 
 of all discrete series representations for $X$.

We address Conjectures \ref{conj:D1} and \ref{conj:D2}
 as subproblems for (A),  
 and formulate Conjecture \ref{conj:D3}
 and Problem \ref{prob:D3}
 for (B).

\begin{problem}
\label{prob:D2}
Find a characterization of the pairs $(G,H)$
 such that $G/H$ admits a discrete series representation.  
\end{problem}

Similarly to the classical fact
 that there is no discrete spectrum
 of the Laplacian $\Delta_{\mathbb{R}^n}$ on ${\mathbb{R}}^n$
 and that there is no continuous spectrum
 of the Laplacian $\Delta_{\mathbb{T}^n}$ 
 on the $n$-torus $\mathbb{T}^n$, 
 the Riemannian symmetric space $G/K$ does not admit 
 any discrete series representation
 if it is of non-compact type
 and does not admit any continuous spectrum
 in the Plancherel formula
 if it is of compact type.  
The answer to Problem \ref{prob:D2}
 is known for reductive symmetric spaces
 by the rank condition:

\begin{equation}
\label{eqn:FJOM}
\operatorname{Disc}(G/H) \ne \emptyset
\iff
\operatorname{rank} G/H = \operatorname{rank} K/H \cap K.  
\end{equation}
The equivalence \eqref{eqn:FJOM} was proved 
 by Flensted-Jensen for $\Leftarrow$
 and Matsuki--Oshima for $\Rightarrow$.  
It generalizes the Riemannian case $G/K$
 as well as Harish-Chandra's rank condition 
 for a group manifold, 
 see \cite{MtO84}
 and references therein.

Beyond symmetric spaces, 
 several approaches
 ({\it{e.g.}}, branching laws \cite{KInvent94, xkdisc}, 
 the wave front set \cite{HOpre}, 
etc.)
 have been applied to find new families
 of
(not necessarily,  real spherical)
 homogeneous spaces
 $G/H$
 that admit discrete series representations.  
It is more involved to prove the opposite direction, 
 {\it{i.e.}}, 
 to prove $\operatorname{Disc}(G/H) = \emptyset$
 for non-symmetric spaces, 
 and very little is known so far, 
 except for certain families of spherical homogeneous spaces.  
For instance, 
 one has:

\begin{example}
[real forms of $SL(2n+1, {\mathbb{C}})/Sp(n, {\mathbb{C}})$, 
{\cite{KInvent94}}]
$\operatorname{Disc}(G/H) = \emptyset$
 if 
$G/H=SL(2n+1, {\mathbb{R}})/Sp(n,{\mathbb{R}})$, 
 whereas 
$\# \operatorname{Disc}(G/H) = \infty$
 for other real forms
 of $G_{\mathbb{C}}/H_{\mathbb{C}}$, 
 {\it{i.e.}}, 
 $SU(2p, 2q+1)/Sp(p,q)$
 or $SU(n,n+1)/Sp(n,{\mathbb{R}})$.  
\end{example}

An optimistic solution to Problem \ref{prob:D2} may be a combination
 of the following two conjectures:

\begin{conjecture}
[{\cite[Conj.\ 6.9]{xkrons94}}]
\label{conj:D1}
In Setting \ref{set:reductive}, 
 one has the equivalence:
\[
\operatorname{Disc}(G/H) \ne \emptyset
\iff
\# \operatorname{Disc}(G/H)=\infty.
\]
\end{conjecture}

\begin{conjecture}
\label{conj:D2}
In Setting \ref{set:reductive}, 
 one has the following equivalence:
\[
\# \operatorname{Disc}(G/H) =\infty
\iff
\text{${\mathfrak{h}}^{\perp} \cap {\mathfrak{g}}_{\operatorname{ell}}^{\ast}$
contains a non-empty open set of ${\mathfrak{h}}^{\perp}$.  }
\]
\end{conjecture}

Both of the conjectures are true
 for reductive symmetric spaces $G/H$.  
In fact, 
 Conjecture \ref{conj:D2} is a reformulation
 of the rank condition \eqref{eqn:FJOM}
 in the spirit of the orbit philosophy.

There are counterexamples
 for the implication $\Rightarrow$
 of an analogous statement
 to Conjectures \ref{conj:D1} and \ref{conj:D2}
 if we drop the assumption 
 that $H$ is reductive,
for instance, they fail when $H$ is a parabolic sugroup and a cocompact discrete subgroup
of $G$ with $\operatorname{rank} G > \operatorname{rank} K$,
respectively. 
The implication $\Leftarrow$ in Conjecture \ref{conj:D2} 
has been proved in  Harris--Y.\ Oshima \cite{HOpre} recently
 without reductivity assumption on $H$.  

\begin{remark}
\label{rem:temp}
Similarly to Conjecture \ref{conj:D2}, 
 one might expect
 the equivalence:
\[
  \text{$L^2(G/H)$ is tempered}
  \iff
  \text{${\mathfrak{h}}^{\perp} \cap {\mathfrak{g}}_{\operatorname{reg}}^{\ast}$ is dense in ${\mathfrak{h}}^{\perp}$.  }
\]
This is proved in Benoist--Kobayashi \cite{BKpre}
 for complex homogeneous spaces
 for any algebraic subgroup $H$
 without reductivity assumption.  
\end{remark}

Once we know $\operatorname{Disc}(G/H) \ne \emptyset$, 
 we may wish to capture {\it{all}} elements
 of $\operatorname{Disc}(G/H)$.  
We divide this exhaustion problem
 into two questions:
one is geometric (Conjecture \ref{conj:D3}), 
and the other is algebraic (Problem \ref{prob:D3}).  

\begin{conjecture}
\label{conj:D3}
Any $\pi \in \operatorname{Disc}(G/H)$ is obtained
 as a geometric quantization
 of some elliptic coadjoint orbit
 that meets ${\mathfrak{h}}^{\perp}$.  
\end{conjecture}

\begin{problem}
\label{prob:D3}
 Find a necessary and sufficient condition
 for cohomologically parabolic induced modules
 $A_{\mathfrak{q}}(\lambda)$ not to vanish
 outside the good range of parameter $\lambda$.  
\end{problem}

Conjecture \ref{conj:D3} strengthens Conjecture \ref{conj:D2}, 
 and one can verify it for reductive symmetric spaces $X$, 
  see \cite[Ex.\ 2.9]{xkrons94}.  
To be more precise, 
 by using 
 Matsuki--Oshima's theorem \cite{MtO84}
 and 
 by using an algebraic characterization
 of Zuckeman's derived functor modules, 
 one can identify any discrete series representation for $G/H$
 as a \lq\lq{geometric quantization}\rq\rq\ $\pi_{\lambda}$
 of an elliptic coadjoint orbit ${\mathcal{O}}_{\lambda}$
 that meets ${\mathfrak{h}}^{\perp}$, 
 with the normalization of \lq\lq{quantization}\rq\rq\
 as in \cite{xkrons94}.  
For \lq\lq{singular}\rq\rq\ $\lambda$, 
 the above $\pi_{\lambda}$ may or may not vanish.  
A missing part of Problem \ref{prob:D1}
 in the literature for symmetric spaces
 is the complete proof of the precise condition
 on $\lambda$
 such that $\pi_{\lambda} \ne 0$, 
 which is reduced 
 to an algebraic question, 
 that is, 
 Problem \ref{prob:D3}.  
The algebraic results in \cite[Chaps.\ 4,5]{K92} and \cite{Tr} 
give an answer
 to Problem \ref{prob:D3} 
 for some classical symmetric spaces.

We examine Problem \ref{prob:D3}
 by $X=X(p,q)$
 with $p,q \ge 1$.  
As we saw in Example \ref{ex:Xpqdisc}, 
 the underlying $({\mathfrak{g}}, K)$-modules (see \cite[Chap.\ 3]{WaI} for instance)
 of the $L^2$-eigenspace $L^2(X)_{\lambda_k}$
 are expressed by $A_{\mathfrak{q}}(k)$.  
Then there are finitely many exceptional parameters
 $k \in {\mathbb{Z}}$
 satisfying $-\frac  12 (p+q-1)<k <0$, 
 {\it{i.e., }}
 lying \lq\lq{outside the good range}\rq\rq\
 for which the general algebraic representation theory
 does not guarantee the irreducibility/non-vanishing for the cohomological parabolic induction.  
This is the point
 that Problem \ref{prob:D3} concerns with.

\section{Problems on discontinuous groups for $G/H$}
\label{sec:discgp}

The local to global study of geometries was a major trend of 20th 
century geometry, with remarkable developments achieved particularly in 
Riemannian geometry. In contrast, in areas such as pseudo-Riemannian 
geometry, familiar to us as the space-time of relativity theory, and 
more generally in manifolds
 with indefinite metric tensor of arbitrary signature, 
surprisingly little is known about global properties of the geometry.  
For instance, 
 the pseudo-Riemannian space form problem is unsolved, 
 which asks the existence of 
 a compact pseudo-Riemannian manifold $M$
 with constant sectional curvature
 for a given signature $(p,q)$, 
 see Conjecture \ref{conj:G4} below.

When we highlight
 \lq\lq{homogeneous structure}\rq\rq\
 as a local property, 
 \lq\lq{discontinuous groups}\rq\rq\
 are responsible
 for the global geometry.  
The theory of discontinuous groups
 beyond the Riemannian setting
 is a relatively \lq\lq{young area}\rq\rq\
 in group theory
 that interacts 
 with topology, 
 differential geometry, 
 representation theory, 
 and number theory among others.  
See \cite{K97} for some background
 on this topic at an early stage
 of the developments.  
This theme was discussed also 
 as a new topic
 of future research
 at the occasion 
 of the World Mathematical Year 2000
 by Kobayashi \cite{K01} and Margulis \cite{Mr}.  
For over 30 years, 
 there have been remarkable developments
 by various methods ranging from topology
 and differential geometry
 to representation theory and ergodic theory, 
 however, 
 some important problems
 are still unsolved, 
 which we illustrate in this section
 by using simple examples.

Beyond the Riemannian setting, 
 we highlight
 a substantial difference
 in \lq\lq{discrete subgroups}\rq\rq\
 and \lq\lq{discontinuous groups}\rq\rq, 
 {\it{e.g.,}} \cite{K89}.

\begin{definition}
\label{def:G1}
Let $G$ be a Lie group acting
 on a manifold $X$.  
A discrete subgroup $\Gamma$ of $G$ 
 is said to be a {\it{discontinuous group}} for $X$
 if $\Gamma$ acts properly discontinuously and freely on $X$.  
\end{definition}

The quotient space $X_{\Gamma}:=\Gamma \backslash X$
 by a discontinuous group $\Gamma$ is 
 a (Hausdorff) $C^{\infty}$-manifold, 
 and any $G$-invariant local geometric structure on $X$
 can be pushed forward to $X_{\Gamma}$
 via the covering map $X \to X_{\Gamma}$.  
Such quotients $X_{\Gamma}$ are complete $(G,X)$-manifolds
 in the sense of Ehresmann and Thurston.

A classical example is a compact Riemann surface $\Sigma_g$
 with genus $g \ge 2$, 
 which can be expressed 
 as $X_{\Gamma}$
 where $\Gamma \simeq \pi_1(\Sigma_g)$
 (surface group)
 and $X \simeq PSL(2,{\mathbb{R}})/PSO(2)$
 by the uniformization theory.  
More generally, 
 any complete affine locally symmetric space
 is given as the form $\Gamma\backslash G/H$
 where $\Gamma$ is a discontinuous group
 for a symmetric space $G/H$.

\begin{remark}
The crucial assumption in Definition \ref{def:G1} is proper discontinuity
 of the action, 
 and freeness is less important.  
In \cite[Def.\ 2.5]{K97}, 
 we did not include the freeness assumption
 in the definition of discontinuous groups, 
 allowing $X_{\Gamma}=\Gamma \backslash X$ to be an orbifold.  
\end{remark}

We discuss the following problems
 in the generality of Setting \ref{set:reductive}, 
 {\it{cf.}} \cite[Problems B and C]{K01}.  

\begin{problem}
\label{prob:G1}
Determine all pairs $(G,H)$ such that $G/H$ admits cocompact discontinuous groups.  
\end{problem}

\begin{problem}
[higher Teichm{\"u}ller theory for $G/H$]
\label{prob:G4}
Describe the moduli of all deformations of a discontinuous group
 $\Gamma$ for $G/H$. 
\end{problem}

In the classical case where $H$ is compact, 
 a theorem of Borel answers Problem \ref{prob:G1} in the affirmative 
 by the existence of cocompact arithmetic discrete subgroups in $G$, 
 whereas the Selberg--Weil local rigidity theorem tells
 that Problem \ref{prob:G4} makes sense for a cocompact $\Gamma$
 in a simple Lie group $G$
 only if ${\mathfrak{g}} \simeq {\mathfrak{s l}}(2,{\mathbb{R}})$, 
 and in this case the deformation
 of discontinuous groups
 gives rise to that of complex structures
 on the Riemann surface.

Such features
 change dramatically
 when $H$ is non-compact:
 some homogeneous spaces may not admit any discontinuous group
 of infinite order 
 (the Calabi--Markus phenomenon \cite{CM}), 
 showing an obstruction to the existence
 of cocompact discontinuous groups for $G/H$.  
On the other hand, 
 discontinuous groups 
 for pseudo-Riemannian manifolds $G/H$
 tend to be \lq\lq{more flexible}\rq\rq\
 in contrast to the classical rigidity theorems 
 in the Riemannian case.  
For instance, 
 some irreducible symmetric spaces
 of arbitrarily higher dimension
 admit cocompact discontinuous groups
 which are not locally rigid
 (\cite{Kas12, K98}), 
 providing wide open settings 
 for Problem \ref{prob:G4}.

As we mentioned, 
 the notion \lq\lq{discontinuous group for $G/H$}\rq\rq\
 is much stronger than 
 \lq\lq{discreteness in $G$}\rq\rq\
 when $H$ is non-compact.  
For instance, 
 a cocompact discrete subgroup $\Gamma$ of $G$ never acts 
 properly discontinuously on $G/H$ unless $H$ is compact.  
Thus the existence of a lattice
 in $G$ does not imply
 that $G/H$ admits a cocompact discontinuous group.

We examine some related questions and conjectures
 to Problem \ref{prob:G1}.  
First, 
 by relaxing the 
 \lq\lq{cocompactness}\rq\rq\ assumption
 of $\Gamma$ in Problem \ref{prob:G1}, 
one may ask the following:

\begin{problem}
\label{prob:G2}
Find a necessary and sufficient condition
 for $G/H$ in Setting \ref{set:reductive}
 to admit a discontinuous group $\Gamma$
 for $G/H$
 such that 
\newline\indent
{\rm{(1)}}\enspace
$\Gamma \simeq {\mathbb{Z}}$;
\newline\indent
{\rm{(2)}}
$\Gamma \simeq \text{a surface group $\pi_1(\Sigma_g)$ with $g \ge 2$.}$  
\end{problem}

Problem \ref{prob:G2} (1) was solved in \cite{K89} 
 in terms of the real rank condition
$\operatorname{rank}_{\mathbb{R}}G >\operatorname{rank}_{\mathbb{R}}H$, 
 which revealed the Calabi--Markus phenomenon \cite{CM}
 in the generality of Setting \ref{set:reductive}.  
Problem \ref{prob:G2} (2) was solved 
 by Okuda \cite{O13}
 for irreducible symmetric spaces, 
 but is unsolved in the generality
 of Setting \ref{set:reductive}.

Cocompact discontinuous groups for $G/H$
 are much smaller
 than cocompact lattices of $G$, 
 for instance, 
 their cohomological dimensions
 are strictly smaller \cite{K89}.  
A simple approach to Problem \ref{prob:G1} is
 to utilize a \lq{continuous analog}\rq\
 of discontinuous groups $\Gamma$:

\begin{definition}
[standard quotients $\Gamma \backslash G/H$ {\cite[Def.\ 1.4]{KasK16}}]
\label{def:standard}
Suppose $L$ is a reductive subgroup of $G$
 such that $L$ acts properly on $G/H$.  
Then any torsion-free discrete subgroup $\Gamma$
 of $L$ is a discontinuous group for $G/H$.  
The quotient space $\Gamma \backslash G/H$
 is called a {\it{standard quotient}} of $G/H$.  
\end{definition}

If such an $L$ acts cocompactly on $G/H$, 
 then $G/H$ admits a cocompact discontinuous group $\Gamma$
 by taking $\Gamma$ to be a torsion-free cocompact discrete subgroup in $L$, 
 which always exists by Borel's theorem.  
We address the following conjecture 
 and a subproblem to Problem \ref{prob:G1}.  

\begin{conjecture}
[{\cite[Conj. 4.3]{K01}}]
\label{conj:G1}
In Setting \ref{set:reductive}, 
 $G/H$ admits a cocompact discontinuous group, 
 only if $G/H$ admits a compact standard quotient.  
\end{conjecture}

If Conjecture \ref{conj:G1} were proved to be true, 
 then Problem \ref{prob:G1} would be reduced
 to the following one:

\begin{problem}
\label{prob:G3}
Classify the pairs $(G,H)$
 such that $G/H$ admits a compact standard quotient.  
\end{problem}

This problem should be manageable because one could use the general theory of real finite-dimensional representations of semisimple Lie algebras and
 apply the properness criterion and the cocompactness criterion in \cite[Thms 4.1 and 4.7]{K89}. See \cite{Tj} for some develepments.

\begin{remark}
(1)\enspace
Conjecture \ref{conj:G1} does not assert
 that any cocompact discontinuous group is a standard one.  
In fact,
 there exist triples $(G, H, \Gamma)$
 such that $\Gamma$ is a cocompact discontinuous group for $G/H$
 and that the Zariski closure of $\Gamma$
 does not act properly on $G/H$, 
 see \cite{Kas12, K98}.  
\par\noindent
(2)\enspace
An analogous statement to Conjecture \ref{conj:G1} fails
 if we drop the reductivity assumption on the groups $G$, $H$ and $L$.  
\par\noindent
(3)\enspace
An analogous statement to Conjecture \ref{conj:G1} is proved
 in Okuda \cite{O13}
 for semisimple symmetric spaces $G/H$
 if we replace the \lq\lq{cocompactness}\rq\rq\ assumption
 with the condition 
 that $\Gamma$ is a surface group $\pi_1(\Sigma_g)$.  
\end{remark}

Special cases of Conjecture \ref{conj:G1} include:
\begin{conjecture}
\label{conj:SLSL}
$SL(n,{\mathbb{R}})/SL(m,{\mathbb{R}})$ does not admit
 a cocompact discontinuous group for any $n>m$.  
\end{conjecture}

\begin{conjecture}
[Space form conjecture {\cite{KY05}}]
\label{conj:G4}
There exists a compact, complete, 
 pseudo-Riemannian manifold
 of signature $(p,q)$
 with constant sectional curvature $1$
 if and only if $(p,q)$ is in the list of Example \ref{ex:GammaXpq} (4) below.  
\end{conjecture}

A criterion on triples $(G,H,L)$
 of reductive Lie groups 
 for $L$ 
 to act properly on $X=G/H$ was established in \cite{K89}, 
 and a list of irreducible symmetric spaces $G/H$
 admitting proper and cocompact actions
 of reductive subgroups $L$
 was given in \cite{KY05}.  
Tojo \cite{Tj} worked with simple Lie groups $G$
 and announced
 that the list in \cite{KY05} exhausts
 all such triples $(G,H,L)$ with $L$ maximal, 
 giving a solution to Problem \ref{prob:G3}
 for symmetric spaces $G/H$
 with $G$ simple.

A number of obstructions 
 to the existence of cocompact discontinuous groups
 for $G/H$
 with $H$ non-compact have been found for over 30 years.  
One of the recent developments includes the affirmative solution
 to the \lq\lq{rank conjecture}\rq\rq\
 raised by the author in 1989:
 it was proved in the case $\operatorname{rank}G=\operatorname{rank} H$ 
 by Kobayashi--Ono (1990), 
 and has been proved recently in the general case by Morita \cite{mor19}
 and Tholozan \cite{Th}, 
 independently.

\begin{conjecture}
[{\cite[Conj.\ 4.15]{K97}}]
\label{conj:G2}
If $G/H$ admits a cocompact discontinuous group, 
 then
$
\operatorname{rank} G + \operatorname{rank} (H \cap K)
\ge 
\operatorname{rank} H + 
\operatorname{rank} K.  
$
\end{conjecture}

Whereas the idea of standard quotients $\Gamma \backslash G/H$
 is to replace a discrete subgroup $\Gamma$ 
with a connected subgroup $L$
 (Definition \ref{def:standard}), 
 one may consider an \lq\lq{approximation}\rq\rq\
 of Problem \ref{prob:G1}, 
 by taking the {\it{tangential homogeneous space}}
 $X_{\theta}=G_{\theta}/H_{\theta}$ 
 in replacement of $X=G/H$, 
 where $G_{\theta} := K \ltimes {\mathfrak{p}}$
 is the Cartan motion group 
 of the real reductive group $G=K \exp {\mathfrak{p}}$
 and similarly for $H_{\theta}$.  
If $G/H$ admits a compact standard quotient, 
 then the tangential homogeneous space
 $G_{\theta}/H_{\theta}$
 admits a cocompact discontinuous group.  
The group $G_{\theta}$ is a compact extension
 of the abelian group ${\mathfrak{p}}$, 
 and is of much simpler structure.  
We ask the following digression
 of Problem \ref{prob:G1}:

\begin{problem}
[{\cite{KY05}}]
\label{prob:G5}
For which pairs $(G,H)$ in Setting \ref{set:reductive}, 
 does $G_{\theta}/H_{\theta}$ admit a cocompact
 discontinuous group?
\end{problem}

This problem is unsolved
 even for symmetric spaces in general, 
 but has a complete answer in some special settings, 
 {\it{e.g.}}, Example \ref{ex:GammaXpq} (6) below.

We end this section
 with a brief summary of the state-of-art
 for these problems
 and conjectures 
 by taking the space form $X(p,q)$ 
 as an example.

\begin{example}
[{\cite{CM, Kan, Kas12, K97, K98, KY05, Ku, mor19, O13, Th}}]
\label{ex:GammaXpq}
Let $(G,H)=(O(p+1,q), O(p,q))$,  
 and $X=X(p,q)=G/H$  the pseudo-Riemannian space form
 of signature $(p,q)$ 
 as in Example \ref{ex:1.4}.  
\par\noindent
(1)\enspace
$X(p,q)$ admits a discontinuous group
 of infinite order iff $p<q$.  
\par\noindent
(2)\enspace
$X(p,q)$ admits a discontinuous group 
 which is isomorphic to a surface group
 iff $p+1 <q$
 or $p+1=q \in 2 {\mathbb{N}}$.  
\par\noindent
(3)\enspace
If $X(p,q)$ admits a cocompact discontinuous group, 
 then $p=0$ or
 \lq\lq{$p<q$ and $q \in 2{\mathbb{N}}$}\rq\rq.  
\par\noindent
(4)\enspace
$X(p,q)$ admits a cocompact discontinuous group
 if $(p,q)$ is in the list below.  
(The converse assertion was stated as Conjecture \ref{conj:G4}.)

\begin{figure}[H]
\begin{center}
\begin{tabular}{cccccc}
\hline
$p$
& ${\mathbb{N}}$
& 0
& 1
& 3
& 7
\\
\hline
$q$
& 0
& ${\mathbb{N}}$
& $2{\mathbb{N}}$
& $4{\mathbb{N}}$
& 8
\\
\hline
\end{tabular}
\end{center}
\end{figure}

\par\noindent
(5)\enspace
If $(p,q)=(0,2)$, $(1,2)$, or $(3,4)$, 
 then $X(p,q)$ admits a cocompact discontinuous group
 which can be deformed continuously
 into a Zariski dense subgroup of $G$
 by keeping proper discontinuity of the action.  
For $(p,q)=(1,2n)$ $(n \ge 2)$, 
 the anti-de Sitter space $X(1,2n)$ admits a compact quotient
 which has a non-trivial continuous deformation
 as standard quotients.  
\par\noindent
(6)\enspace
The tangential homogeneous space
 $G_{\theta}/H_{\theta}$ admits a cocompact discontinuous group
 if and only if 
 $p< \rho(q)$ 
 where $\rho(q)$ is the Radon--Hurwitz number, 
 or equivalently, 
 if and only if $(p,q)$ is in the following list:

\begin{figure}[H]
\begin{center}
\begin{tabular}{cccccccccccccc}
\hline
$p$
& ${\mathbb{N}}$
& 0
& 1
& 3
& 4
& 5
& 6
& 7
& 8
& 9 
& 10
& 11
& $\cdots$
\\
\hline
$q$
& 0
& ${\mathbb{N}}$
& $2{\mathbb{N}}$
& $4{\mathbb{N}}$
& $8{\mathbb{N}}$
& $8{\mathbb{N}}$
& $8{\mathbb{N}}$
& $8{\mathbb{N}}$
& $16{\mathbb{N}}$
& $32{\mathbb{N}}$
& $64{\mathbb{N}}$
& $64{\mathbb{N}}$
& $\cdots$
\\
\hline
\end{tabular}
\end{center}
\end{figure}

\end{example}

\section{Spectral analysis for pseudo-Riemannian locally homogeneous spaces
 $\Gamma \backslash G/H$}
\label{sec:spec}

This section discusses briefly a new direction of analysis
 on pseudo-Riemannian locally homogeneous spaces $\Gamma \backslash G/H$.

Suppose we are in Setting \ref{set:reductive}.  
Let $\Gamma$ be a discontinuous group for $X=G/H$.  
Then any $G$-invariant differential operator
 $D \in {\mathbb{D}}(G/H)$ induces 
 a differential operator $D_{\Gamma}$
 on the quotient $X_{\Gamma}:=\Gamma \backslash G/H$
 via the covering $X \to X_{\Gamma}$.  
We think of the set
 ${\mathbb{D}}(X_{\Gamma}):=\{D_{\Gamma}: D \in  {\mathbb{D}}(G/H)\}$
 as the algebra
 of {\it{intrinsic differential operators}}
 on the locally homogeneous space $X_{\Gamma}$.  

\begin{example}
\label{ex:Lap}
(1)\enspace
In Setting \ref{set:reductive}, 
 $X_{\Gamma}$ inherits a pseudo-Riemannian structure from $X$, 
 and the Laplacian $\Delta_{X_{\Gamma}}$
 belongs to ${\mathbb{D}}(X_{\Gamma})$.  
\par\noindent
(2)\enspace
For $X=X(p,q)$, 
 ${\mathbb{D}}(X_{\Gamma})$ is a polynomial ring in
 the Laplacian $\Delta_{X_{\Gamma}}$
 for any discontinuous group $\Gamma$.  
\end{example}

We address the following problem:

\begin{problem}
[see {\cite{KasK16, KasKpre}}]
\label{prob:S1}
For intrinsic differential operators on $X_{\Gamma}=\Gamma \backslash G/H$, 
\par\noindent
(1)\enspace
construct joint eigenfunctions on $X_{\Gamma}$;
\par\noindent
(2)\enspace
find a spectral theory on $L^2(X_{\Gamma})$.  
\end{problem}

With the same spirit as in Section \ref{sec:discrep}, 
 we highlight \lq\lq{discrete spectrum}\rq\rq.   

\begin{definition}
We say $\lambda \in \invHom{{\mathbb{C}}\text{-alg}}{{\mathbb{D}}(X_{\Gamma})}{\mathbb{C}}$ is a {\it{discrete spectrum}}
 for intrinsic differential operators on $X_{\Gamma}$
 if $L^2(X_{\Gamma})_{\lambda} \ne \{0\}$, 
 where we set
\[
   L^2(X_{\Gamma})_{\lambda}
   :=
   \{f \in L^2(X_{\Gamma}): D f = {\lambda}(D) f
\quad 
{}^{\forall} D \in {\mathbb{D}}(X_{\Gamma})\}.  
\]
We write $\operatorname{Spec}_d(X_{\Gamma})$
 for the set of discrete spectra.  
\end{definition}

A subproblem to Problem \ref{prob:S1} (1) includes:

\begin{problem}
\label{prob:S2}
Construct joint $L^2$-eigenfunctions on $X_{\Gamma}$.  
\end{problem}

In relation to Problem \ref{prob:G4} about the deformations
 of a discontinuous group $\Gamma$ for $G/H$, 
 one may also ask the following:

\begin{problem}
\label{prob:S3}
Understand the behavior of $\operatorname{Spec}_d(X_{\Gamma})$ 
 under small deformations of $\Gamma$ inside $G$.  
\end{problem}

These problems have been 
studied extensively
 in the following special settings
 for $X_{\Gamma}=\Gamma\backslash G/H$:

\begin{enumerate}
\item[(1)]
$(H=K)$.  
When $H$ is a maximal compact subgroup $K$ of $G$, 
 {\it{i.e.,}}
 $X_{\Gamma}$ is a {\it{Riemannian}} locally symmetric space,
 a vast theory has been developed
 over several decades, 
 in particular, 
 in connection with the theory 
 of automorphic forms
 when $\Gamma$ is arithmetic.  
\item[(2)]
($\Gamma=\{e\}$).  
This case is related to the topic in Section \ref{sec:discrep}.  
In particular, 
 Problem \ref{prob:S1} has been extensively studied in the case
 where $G/H$ is a reductive symmetric space
 and $\Gamma=\{e\}$.  
\item[(3)]
$G={\mathbb{R}}^{p,q}$, 
 $\Gamma={\mathbb{Z}}^{p+q}$, 
 and $H=\{0\}$.  
In this case, 
 $\operatorname{Spec}_d(X_{\Gamma})$ is the set
 of values of indefinite quadratic forms
 at integral points, 
 see \cite{K16}
 for a discussion on Problem \ref{prob:S3}
 in relation to the Oppenheim conjecture
 (proved by Margulis)
 in Diophantine approximation.  
\end{enumerate}

The situation changes drastically
 beyond the classical setting, 
 namely, 
 when $H$ is not compact any more and $\Gamma \ne \{e\}$.  
New difficulties include:

\begin{enumerate}
\item[(1)]
(Representation theory)\enspace
Even when $\Gamma \backslash G/H$ is compact, 
 the regular representation of $G$ 
 on $L^2(\Gamma \backslash G)$
 has infinite multiplicities, 
 as opposed to a classical theorem 
 of Gelfand--Piatetski--Shapiro.  
\item[(2)]
(Analysis)\enspace
In contrast to the Riemannian case
 where $H=K$, 
 the Laplacian $\Delta_{X_{\Gamma}}$ is not an elliptic 
differential operator any more.  
\end{enumerate}

As we saw in Section \ref{sec:discgp}, 
 if $H$ is not compact, 
 then not all homogeneous spaces $G/H$ admit discontinuous groups
 of infinite order, 
 but fortunately, 
 there exist a family of reductive symmetric spaces $G/H$
 that admit \lq\lq{large}\rq\rq\ discontinuous groups $\Gamma$, 
 {\it{e.g.,}} such that $X_{\Gamma}=\Gamma \backslash G/H$ is compact 
or of finite volume.  
Moreover, 
there also exist triples $(G,H,\Gamma)$
such that discontinuous groups $\Gamma$ for $G/H$
 can be deformed continuously.  
These examples offer broad settings
 for Problems \ref{prob:S1}
 and its subproblems.

For Problem \ref{prob:S3}, 
 we consider two notions for stability:

\begin{definition}
(1) (stability for proper discontinuity)\enspace
A discontinuous group $\Gamma$
 is {\it{stable under small deformations}}
 if the group $\varphi(\Gamma)$ acts properly discontinuously
 and freely on $X$
 for all $\varphi \in \invHom{}{\Gamma}{G}$
 in some neighbourhood ${\mathcal{U}}$
 of the natural inclusion $\Gamma$ in $G$.  
\par\noindent
(2) (stability for $L^2$-spectrum)\enspace
We say $\lambda \in \invHom{{\mathbb{C}}\text{-alg}}{{\mathbb{D}}(X_{\Gamma})}{\mathbb{C}}$
 is a {\it{stable spectrum}}
 if $L^2(X_{\varphi(\Gamma)})_{\lambda} \ne \{0\}$
 for any $\varphi \in \invHom{}{\Gamma}{G}$
 in some neighbourhood ${\mathcal{U}}$
 of the natural inclusion $\Gamma$ in $G$.  
\end{definition}

\begin{conjecture}
\label{conj:S1}
Suppose that $\Gamma$ is a finitely generated discontinuous group for $G/H$
 having
 non-trivial continuous deformations (up to inner automorphisms) 
 with stability of proper discontinuity.  
Then the following conditions on the pair $(G,H)$ 
 are equivalent.  
\par\noindent
(i)\enspace
There exist infinitely many stable spectra on $L^2(\Gamma \backslash G/H)$.  

\par\noindent
(ii)\enspace
$\operatorname{Disc}(G/H) \ne \emptyset$.  
\end{conjecture}

See \cite{KasK16}
 for some results
 in the direction (ii) $\Rightarrow$ (i), 
 which treat also the case
 $\operatorname{vol}(\Gamma \backslash G/H)=\infty$.

The last section has been devoted
 to a \lq\lq{very young}\rq\rq\ topic, 
 though special cases trace back to rich and deep classical theories.  
I expect this topic will create new interactions
 with different subjects of mathematics, 
 and this is why I include it as a part
 of my article for \lq\lq{Mathematics Going Forward}\rq\rq.

\vskip 1.5pc
\par\noindent
{\bf{$\langle$Acknowledgements$\rangle$}}\enspace
The author would like to express
 his sincere gratitude 
 to his collaborators
 of the various projects mentioned in this article.  
This work was partially supported
 by Grant-in-Aid for Scientific Research (A) (18H03669), 
JSPS.

\end{document}